\documentclass[12pt,a4paper]{article}
\usepackage{amsmath}
\usepackage{amssymb}
\usepackage{t1enc}
\usepackage[latin1]{inputenc}
\usepackage[english]{babel}

\usepackage{mathtools}

\usepackage{amsfonts}
\usepackage{latexsym}
\usepackage{bm}
\newtheorem{theorem}{Theorem}[section]

\newtheorem{lemma}[theorem]{Lemma}

\allowdisplaybreaks

\begin{document}
\title{Reducing the maximum degree of a graph: comparisons of bounds}

\author{Peter Borg\\[5mm]
Department of Mathematics \\
Faculty of Science \\
University of Malta\\
Malta\\
\texttt{peter.borg@um.edu.mt} 
}

\date{}
\maketitle

\begin{abstract}
Let $\lambda(G)$ be the smallest number of vertices that can be removed from a non-empty graph $G$ so that the resulting graph has a smaller maximum degree. Let $\lambda_{\rm e}(G)$ be the smallest number of edges that can be removed from $G$ for the same purpose. Let $k$ be the maximum degree of $G$, let $t$ be the number of vertices of degree $k$, let $M(G)$ be the set of vertices of degree $k$, let $n$ be the number of vertices in the closed neighbourhood of $M(G)$, and let $m$ be the number of edges incident to vertices in $M(G)$. Fenech and the author showed that $\lambda(G) \leq \frac{n+(k-1)t}{2k}$, and they essentially showed that $\lambda (G) \leq n \left ( 1- \frac{k}{k+1} { \Big( \frac{n}{(k+1)t} \Big) }^{1/k} \right )$. They also showed that $\lambda_{\rm e}(G) \leq \frac{m + (k-1)t}{2k-1}$ and $\lambda_{\rm e} (G) \leq m \left ( 1- \frac{k-1}{k} { \Big( \frac{m}{kt} \Big) }^{1/(k-1)} \right )$. These bounds are attained if $k \geq 2$ and $G$ is the union of $t$ pairwise vertex-disjoint $(k+1)$-vertex stars. For each of $\lambda(G)$ and $\lambda_{\rm e}(G)$, the two bounds on the parameter are compared for the purpose of determining, for each bound, the cases in which the bound is better than the other. This work is also motivated by the likelihood that similar pairs of bounds will be discovered for other graph parameters and the same analysis can be applied. 
\end{abstract}

\section{Introduction} \label{introsection}

For basic terminology and notation in graph theory, we refer the reader to \cite{Bollobas, West}; the definitions of terms and notations used here are given in the papers \cite{BF1,BF2}, which are the basis of the work presented here. 

The set $\{1, 2, \dots\}$ of positive integers is denoted by $\mathbb{N}$. For any $n \in \mathbb{N}$, the set $\{1, \dots, n\}$ is denoted by $[n]$. For a set $X$, the set of all $2$-element subsets of $X$ is denoted by ${X \choose 2}$. Arbitrary sets are taken to be finite. 

Every graph $G$ is taken to be \emph{simple}, that is, its vertex $V(G)$ and edge set $E(G)$ satisfy $E(G) \subseteq {V(G) \choose 2}$. We may represent an edge $\{v,w\}$ by $vw$. For $v \in V(G)$, $N_{G}(v)$ denotes $\{w \in V(G) \colon vw \in E(G)\}$, $N_{G}[v]$ denotes $N_{G} (v) \cup \{ v \}$,  $E_G(v)$ denotes $\{e \in E(G) \colon v \in e\}$, and $d_{G} (v)$ denotes $|N_{G} (v)|$ ($= |E_G(v)|$) and is called the \emph{degree of $v$}. For $X \subseteq V(G)$, $\bigcup_{v \in X} N_G[v]$ is denoted by $N_G[X]$ and called the \emph{closed neighbourhood of $X$}. The \emph{maximum degree of $G$} is $\max \{d_{G}(v) \colon v \in V(G)\}$ and is denoted by $\Delta(G)$. The set of vertices of $G$ of degree $\Delta(G)$ is denoted by $M(G)$. For $X \subseteq V(G)$, $G[X]$ denotes the \emph{subgraph of $G$ induced by $X$}, that is, $G[X] = (X,E(G) \cap {X \choose 2})$. For $R \subseteq V(G)$, $G-R$ denotes the subgraph of $G$ obtained by removing the vertices in $R$ from $G$, that is, $G-R = G[V(G) \backslash R]$. For $L \subseteq E(G)$, $G - L$ denotes the subgraph of $G$ obtained by removing the edges in $L$ from $G$, that is, $G - L = (V(G), E(G) \backslash L)$. 

We call a subset $R$ of $V(G)$ a \emph{$\Delta$-reducing set of $G$} if $\Delta(G-R) < \Delta(G)$ or $R = V(G)$ (note that $V(G)$ is the smallest $\Delta$-reducing set of $G$ if and only if $E(G) = \emptyset$). 
We call a subset $L$ of $E(G)$ a \emph{$\Delta$-reducing edge set of $G$} if $\Delta(G-L) < \Delta(G)$ or $L = E(G) = \emptyset$. We denote the size of a smallest $\Delta$-reducing set of $G$ by $\lambda(G)$, and we denote the size of a smallest $\Delta$-reducing edge set of $G$ by $\lambda_{\rm e}(G)$.

Let $G_{\rm v}$ denote the subgraph of $G$ induced by $N_G[M(G)]$, and let $G_{\rm e}$ denote the subgraph of $G$ with vertex set $N_G[M(G)]$ and edge set $\bigcup_{v \in M(G)} E_G(v)$. As explained in \cite{BF1, BF2}, we clearly have
\begin{gather} \Delta(G_{\rm v}) = \Delta(G), \quad M(G_{\rm v}) = M(G), \quad \lambda(G_{\rm v}) = \lambda(G), \label{eqnsve1} \\
\Delta(G_{\rm e}) = \Delta(G), \quad M(G_{\rm e}) = M(G), \quad \lambda_{\rm e}(G_{\rm e}) = \lambda_{\rm e}(G), \label{eqnsve2} \\
|V(G_{\rm v})| \leq \sum_{v \in M(G)} |N_G[v]| = (\Delta(G) +1)|M(G)|, \label{eqnsve3} \\
|E(G_{\rm e})| \leq \sum_{v \in M(G)} |E_G(v)| = \Delta(G)|M(G)|. \label{eqnsve4}
\end{gather}
By the handshaking lemma and (\ref{eqnsve2}), $2|E(G_{\rm e})| = \sum_{v \in V(G_{\rm e})} d_{G_{\rm e}}(v) \geq \sum_{v \in M(G_{\rm e})} d_{G_{\rm e}}(v) = \Delta(G_{\rm e})|M(G_{\rm e})| = \Delta(G)|M(G)|$, so
\begin{equation} |E(G_{\rm e})| \geq \Delta(G)|M(G)|/2. \label{eqnsve5}
\end{equation}

The graph parameters $\lambda(G)$ and $\lambda_{\rm e}(G)$ were investigated in \cite{BF1} and \cite{BF2}, respectively. For each of them, two main general bounds were obtained, and the bounds are sharp. 

The following is the first main general bound proved in \cite{BF1}.

\begin{theorem}[\cite{BF1}] \label{v1}
If $G$ is a graph, $k = \Delta(G) \geq 1$, $t = |M(G)|$, and $n = |V(G_{\rm v})|$, then 
\[\lambda(G) \leq \frac{n+(k-1)t}{2k}.\]
\end{theorem}
In \cite{BF1}, the result is actually stated with $n = |V(G)|$, but the improvement given by $n = |V(G_{\rm v})|$ is immediately deduced from (\ref{eqnsve1}). The extremal structures are determined in \cite{BF1b}. 
Using a probabilistic argument similar to that used by Alon in \cite{Alon}, it was also shown in \cite[Proof of Theorem~2.7]{BF1} that 
\begin{equation} \lambda(G) \leq u(p) = np + t(1-p)^{k + 1} \mbox{ for any real number $p$ such that $0 \leq p \leq 1$,} \label{e2ineq}
\end{equation}
and that this yields the bound 
\[\lambda (G) \leq \frac{n\ln{(k +1)} + t}{k+1}.\]
However, by differentiating $u$ with respect to $p$, we find that the minimum value of $u$ occurs at $p = 1 - ( \frac{n}{(k+1)t} )^{1/k}$ (note that this satisfies $0 \leq p \leq 1$ by (\ref{eqnsve3})),  
and hence it is $n(1 - ( \frac{n}{(k+1)t} )^{1/k}) + t( \frac{n}{(k+1)t} )^{1+1/k} = n(1 - \frac{k}{k+1}( \frac{n}{(k+1)t} )^{1/k} )$. Thus, by (\ref{e2ineq}), the following was essentially established in \cite{BF1}.

\begin{theorem} \label{v2}
If $G$ is a graph, $k = \Delta(G) \geq 1$, $t = |M(G)|$, and $n = |V(G_{\rm v})|$, then 
\[\lambda(G) \leq n \left( 1 - \frac{k}{k+1} \left( \frac{n}{(k+1)t} \right)^{1/k} \right).\]
\end{theorem}

The following are the two main general bounds proved in \cite{BF2}.

\begin{theorem}[\cite{BF2}] \label{e1}
If $G$ is a graph, $k = \Delta(G) \geq 1$, $t=|M(G)|$, and $m=|E(G_{\rm e})|$, then
\[\lambda_{\rm e}(G) \leq \frac{m + (k-1)t}{2k-1}. \]
\end{theorem}

\begin{theorem}[\cite{BF2}] \label{e2}
If $G$ is a graph, $k = \Delta(G) \geq 2$, $t=|M(G)|$, and $m=|E(G_{\rm e})|$, then 
\[\lambda_{\rm e} (G) \leq m \left ( 1- \frac{k-1}{k} { \Big( \frac{m}{kt} \Big) }^{1/(k-1)} \right ).\]
\end{theorem}
Theorem~\ref{e2} was obtained by means of a probabilistic argument similar to that for Theorem~\ref{v2}.

The bounds in Theorems~\ref{v1}--\ref{e2} are attained if, for example, $k \geq 2$ and $G$ is the union of $t$ pairwise vertex-disjoint $(k+1)$-vertex stars (see \cite{BF1, BF1b, BF2}). If $k = 1$, then $G_{\rm v}$ and $G_{\rm e}$ are the same union of $t$ pairwise vertex-disjoint $2$-vertex stars, and hence the bounds in Theorems~\ref{v1} and \ref{e1} are attained.

In this paper, we compare the bounds on $\lambda(G)$ in Theorems~\ref{v1} and \ref{v2}, and we compare the bounds on $\lambda_{\rm e}(G)$ in Theorems~\ref{e1} and \ref{e2}. We use several well-known results from real analysis to determine, for each bound, a significant proportion of the cases in which the bound is better than the other bound for the same parameter. Our main contribution is the solution of the problem for $\lambda_{\rm e}(G)$ for $k$ sufficiently large.

\begin{theorem} Let $k \geq 2$, $t$, and $m$ be as in Theorems~\ref{e1} and \ref{e2}, and let $x_1 > 1$ be the real number $3.512...$ such that $e^{(x_1 - 1)/2} = x_1$. \\
(a) If $m = kt$, then the bounds in Theorems~\ref{e1} and \ref{e2} are equal. \\
(b) If $kt/x_1 < m \neq kt$ and $k$ is sufficiently large, then the bound in Theorem~\ref{e1} is smaller than the bound in Theorem~\ref{e2}. \\
(c) If $m <  kt/x_1$ and $k$ is sufficiently large, then the bound in Theorem~\ref{e2} is smaller than the bound in Theorem~\ref{e1}.
\end{theorem}
%
A stronger version that addresses any $k \geq 2$ is proved in Section~\ref{comp1sectione}.

This work is also motivated by the likelihood that similar pairs of bounds will be discovered for other graph parameters and the same analysis can be applied. 

\section{The bounds in Theorems~\ref{e1} and \ref{e2}} \label{comp1sectione}

Let $k \geq 2$, $t$, and $m$ be as in Theorems~\ref{e1} and \ref{e2}. Let $b_1(k,t,m)$ and $b_2(k,t,m)$ be the bound in Theorem~\ref{e1} and the bound in Theorem~\ref{e2}, respectively; that is,
\[b_1(k,t,m) = \frac{m + (k-1)t}{2k-1} \quad \mbox{and} \quad b_2(k,t,m) = m \left ( 1- \frac{k-1}{k} { \Big( \frac{m}{kt} \Big) }^{1/(k-1)} \right ).\]
%
By (\ref{eqnsve4}), $m \leq kt$, and equality holds if $G$ is the union of $t$ pairwise vertex-disjoint $(k+1)$-vertex stars, in which case the two bounds are equal and attained. We now consider $m < kt$. 


If $k = 2$, then 
\begin{equation} b_2(k,t,m) - b_1(k,t,m) = m - \frac{m^2}{4t} - \frac{m+t}{3} = \frac{(2t-m)(3m-2t)}{12t} > 0 \nonumber
\end{equation}
as $m < kt = 2t$ and $m \geq kt/2 = t$ by (\ref{eqnsve5}). Thus, $b_1(k,t,m) < b_2(k,t,m)$ if $k = 2$. We now consider $k \geq 3$.

\begin{theorem} \label{compresulte} Suppose $k \geq 3$ and $m < kt$. Let $x_1 > 1$ be the real number $3.512...$ such that $e^{(x_1 - 1)/2} = x_1$.\medskip 
\\
(a) There exists a unique real number $x_0$ such that $2.088... \leq x_0 < x_1$ and $\left( \frac{2k-1}{2k - x_0} \right)^{k} = x_0$. We have 
\[b_1(k,t,m) < b_2(k,t,m) \quad \mbox{if } \; m \geq \frac{1}{x_0} kt.\]
The larger $k$ is, the larger $x_0$ is. Moreover, for any real $\delta > 0$, $x_0 > x_1 - \delta$ if $k$ is sufficiently large.
\medskip 
\\
(b) There exists a unique real number $x_0'$ such that $x_1 \leq x_0' < 4$ and $\left( \frac{2k-2}{2k - 1} \right)^{1/2} e^{(x_0' - 1)/2} = x_0'$. We have 
%
%
\[b_2(k,t,m) < b_1(k,t,m) \quad \mbox{if } \; m \leq \frac{1}{x_0'} kt.\]
The larger $k$ is, the smaller $x_0'$ is. Moreover, for any real $\delta > 0$, $x_0' < x_1 + \delta$ if $k$ is sufficiently large.
\end{theorem}
Since $m$ can be at most $kt$, this result tells us that the range of values of $m$ for which the bound in Theorem~\ref{e1} is better than the bound in Theorem~\ref{e2} is wider than that for which the opposite holds. 

We now prove Proposition~\ref{compresulte}. The set of real numbers is denoted by $\mathbb{R}$, and the set of positive real numbers is denoted by $\mathbb{R}^+$. We shall make use of standard notation for real intervals. Let $e$ be the base of the natural logarithm, that is, $e = \lim_{n \rightarrow \infty} \left( 1 + \frac{1}{n}\right)^n = 2.718...$.

\begin{lemma} \label{limitlemma} If $f: \mathbb{R}^+ \rightarrow \mathbb{R}^+$ is the function given by \[f(x) = \left( 1 + \frac{1}{x} \right)^{x+1}\] for $x > 0$, then $f(x)$ decreases as $x$ increases, and $\lim_{x \rightarrow \infty} f(x) = e$.
\end{lemma}
\textbf{Proof.} Let $g : (-\frac{1}{2}, \infty) \rightarrow \mathbb{R}$ be the function given by \[g(z) = z - \ln(1 + z)\] for $z > -\frac{1}{2}$. The derivative $\frac{{\rm d}g}{{\rm d}z}$ is $1 - \frac{1}{1+z}$, which is negative for $-\frac{1}{2} < z < 0$, $0$ for $z = 0$, and positive for $z > 0$. Thus, $g(z)$ increases from $g(0) = 0$ as $z$ increases from $0$ to infinity, and hence
\begin{equation} g(z) > 0 \quad \mbox{for $z > 0$}. \label{fact1}
\end{equation}

We have $\ln f(x) = (x+1)\ln\left( 1 + \frac{1}{x} \right)$. Using implicit differentiation, we obtain $\frac{1}{f(x)} \frac{{\rm d}f}{{\rm d}x} = \ln\left( 1 + \frac{1}{x} \right) + (x+1)\left(\frac{1}{1 + \frac{1}{x}}\right)\left(-\frac{1}{x^2}\right) = \ln\left( 1 + \frac{1}{x} \right) - \frac{1}{x}$. Thus, by (\ref{fact1}) with $z = \frac{1}{x}$, $-\frac{1}{f(x)} \frac{{\rm d}f}{{\rm d}x} > 0$, and hence, since $f(x) > 0$, we obtain $\frac{{\rm d}f}{{\rm d}x} < 0$. Therefore, $f(x)$ decreases as $x$ increases. Now $\lim_{x \rightarrow \infty} f(x) = \left( \lim_{x \rightarrow \infty} \left( 1 + \frac{1}{x} \right)^{x} \right) \left(\lim_{x \rightarrow \infty} \left( 1 + \frac{1}{x} \right) \right) = e$.~\hfill{$\Box$}

\begin{lemma} \label{limitlemma3} For any $c \in (0,1]$, let $f_c : [1,\infty) \rightarrow \mathbb{R}$ be the function given by \[f_c(x) =  \frac{ce^{(x-1)/2}}{x}\] for $x \geq 1$.
\\
(a) As $x$ increases from $2$ to infinity, $f_c(x)$ increases to infinity.\\
(b) There exists a unique real number $x_c > 1$ such that $f_c(x_c) = 1$, and $x_c > 2$. \\
(c) If $c_1, c_2 \in (0,1]$ with $c_1 < c_2$, then $x_{c_2} < x_{c_1}$. \\
(d) For any real $\delta > 0$, $f_1(x_1 + \delta) > 1$ and $x_c < x_1 + \delta$ for any $c \in (1/f_1(x_1 + \delta), 1]$.
\end{lemma}
\textbf{Proof.} Using differentiation, we obtain that the minimum value of $f_c$ occurs at $x = 2$, and that $f_c$ has no other turning points. Thus, $f_c(x)$ decreases from $f_c(1) = c$ to $f_c(2) = ce^{1/2}/2 < c$ as $x$ increases from $1$ to $2$, and $f_c(x)$ increases as $x$ increases from $2$. Since $f_c(x) = \frac{c}{x}\sum_{i=0}^\infty \frac{((x-1)/2)^i}{i!} = c \left( \frac{5}{8x} + \frac{1}{4} + \frac{x}{8} + \sum_{i=3}^\infty \frac{((x-1)/2)^i}{i!x} \right)$, $f_c(x)$ increases from $f_c(2) < 1$ to infinity as $x$ increases from $2$ to infinity. This yields (a) and (b). 

Let $c_1, c_2 \in (0,1]$ with $c_1 < c_2$. Since $f_{c_2}(x_{c_1}) = \frac{c_2}{c_1} f_{c_1}(x_{c_1}) = \frac{c_2}{c_1} > 1$, we obtain (c). 

Let $\delta \in \mathbb{R}^+$. Let $x' = x_1 + \delta$. Let $y = f_1(x')$. By (a) and (b), $y > 1$. We have $f_{1/y}(x') = \frac{1}{y}f_1(x') = 1$, so $x' = x_{1/y}$. By (c), (d) follows.~\hfill{$\Box$}

\begin{lemma} \label{limitlemma2e} Let $A = \{(x, y) \in \mathbb{R} \times \mathbb{R} \colon y \geq 3, \, 1 \leq x < 2y\}$. Let $f: A \rightarrow \mathbb{R}$ be the function given by
\[f(x,y) = \left( \frac{2y-1}{2y - x} \right)^{y-1/2} - x\] 
for $(x,y) \in A$. For any $y_0 \in [3, \infty)$, $f(x_{y_0}, y_0) = 0$ for some unique $x_{y_0} \in (1, 2y_0)$, and $f(x,y) < 0$ for any $x \in (1,x_{y_0}]$ and $y \in [y_0, \infty)$ such that $x \neq x_{y_0}$ or $y \neq y_0$. 

Moreover, let $x_1 > 1$ be the real number $3.512...$ such that $e^{(x_1 - 1)/2} = x_1$. \\
(a) If $y_0, y_1 \in [3, \infty)$ with $y_0 < y_1$, then $x_{y_0} < x_{y_1} < x_1$. \\
(b) For any real $\delta > 0$, there exists some $y_{\delta} \in [3, \infty)$ such that $x_y > x_1 - \delta$ for any $y \in (y_{\delta}, \infty)$.
\end{lemma}
\textbf{Proof.} Let $g : [1, 2y_0) \rightarrow \mathbb{R}$ such that $g(x) = f(x,y_0)$ for $x \in [1, 2y_0)$. We have
\begin{equation} \frac{{\rm d}g}{{\rm d}x} = (y_0 -1/2) \left( \frac{2y_0 -1}{2y_0 -x} \right)^{y_0 -3/2}\frac{2y_0 -1}{(2y_0 -x)^2}  - 1 = \frac{1}{2}\left( \frac{2y_0 -1}{2y_0 -x} \right)^{y_0 +1/2} - 1. \nonumber
\end{equation} 
As $x$ increases from $1$ to $2y_0$, the value of $\frac{1}{2}\left( \frac{2y_0 -1}{2y_0 -x} \right)^{y_0 +1/2}$ increases from $\frac{1}{2}$ to $\infty$, and hence $\frac{{\rm d}g}{{\rm d}x}$ increases from $-\frac{1}{2}$ to $\infty$. Thus, there exists a unique $x^* \in (1, 2y_0)$ such that $\frac{{\rm d}g}{{\rm d}x}$ is $0$ at $x^*$, and $g(x^*) = \min\{g(x) \colon x \in [1, 2y_0)\} < g(1) = 0$. Thus, $g(x)$ decreases from $g(1) = 0$ to $g(x^*)$, and then increases from $g(x^*)$ to $\infty$. Consequently, there exists a unique $x_{y_0} \in (1, 2y_0)$ such that $g(x_{y_0}) = 0 = g(1)$ and $g(x) < g(x_{y_0})$ for each $x \in (1, x_{y_0})$. 

Now suppose $x \in (1,x_{y_0}]$ and $y \in [y_0, \infty)$. Let $z_0 = \frac{2y_0 - x}{x-1}$ and $z = \frac{2y - x}{x-1}$. Then, $z \geq z_0$. We have
\begin{align} f(x,y) + x &= \left( 1 + \frac{x-1}{2y - x} \right)^{y-1/2} = \left(1 + \frac{1}{z} \right)^{(z+1)(x-1)/2} \nonumber \\
&= \left( \left(1 + \frac{1}{z} \right)^{z+1}\right)^{(x-1)/2} \leq \left( \left(1 + \frac{1}{z_0} \right)^{z_0 +1}\right)^{(x-1)/2} \quad \mbox{(by Lemma~\ref{limitlemma})} \nonumber \\
&= f(x,y_0) + x. \label{fact1.1e}
\end{align} 
Therefore, 
\begin{equation} f(x,y) \leq f(x,y_0) = g(x) \leq g(x_{y_0}) = 0. \label{fact1.2e}
\end{equation}
If $x \neq x_{y_0}$, then $x < x_{y_0}$, and hence $g(x) < g(x_{y_0})$. If $y \neq y_0$, then $y > y_0$, $z > z_0$, $\left(1 + \frac{1}{z} \right)^{z+1} < \left(1 + \frac{1}{z_0} \right)^{z_0 +1}$ (by Lemma~\ref{limitlemma}), and hence $f(x,y) < f(x,y_0)$ by (\ref{fact1.1e}). Thus, if $x \neq x_{y_0}$ or $y \neq y_0$, then $g(x) < g(x_{y_0})$ or $f(x,y) < f(x,y_0)$, and hence $f(x,y) < 0$ by~(\ref{fact1.2e}).

Let $h : [1,\infty) \rightarrow \mathbb{R}$ such that $h(x) = e^{(x-1)/2} - x$ for $x \geq 1$. Using differentiation, we obtain that the minimum value of $h$ occurs at $x = 1 + 2 \ln 2 < x_1$, and that $h$ has no other turning points. Thus, $h(x)$ decreases from $h(1) = 0$ to $h(1 + 2 \ln 2) < 0$ as $x$ increases from $1$ to $1 + 2 \ln 2$, and, since $h(x) = x \left( \frac{e^{(x-1)/2}}{x} - 1 \right)$, Lemma~\ref{limitlemma3} (a) implies that $h(x)$ increases to infinity as $x$ increases from $1 + 2 \ln 2$. Note that $h(x_1) = 0$. Let $z_0' = \frac{2y_0 - x_{y_0}}{x_{y_0}-1}$. We have $0 = f(x_{y_0},y_0) = \left( \left(1 + \frac{1}{z_0'} \right)^{z_0' +1}\right)^{(x_{y_0}-1)/2} - x_{y_0} > e^{(x_{y_0}-1)/2} - x_{y_0}$ by Lemma~\ref{limitlemma}. Thus, $h(x_{y_0}) < 0$, and hence $x_{y_0} < x_1$.

Next, suppose $y_0 < y_1$. By the same argument for $x_{y_0}$, $x_{y_1} < x_1$. Let $p : [1, 2y_1) \rightarrow \mathbb{R}$ such that $p(x) = f(x,y_1)$ for $x \in [1, 2y_1)$. By the argument above, $p(x) < 0$ for $x \in (1,x_{y_1})$, $p(x_{y_1}) = 0$, and $p(x) > 0$ for $x \in (x_{y_1}, 2y_1)$. 
Since $y_1 \in (y_0, \infty)$, we have $p(x) = f(x,y_1) < 0$ for any $x \in (1,x_{y_0}]$, so $x_{y_1} > x_{y_0}$. Thus, (a) is proved.

Finally, let $\delta \in \mathbb{R}^+$. Let $\delta' = \min\{\frac{1}{2}, \delta\}$. Let $x' = x_1 - \delta'$. Since $x' < x_1$, $h(x') < 0$. Let $q : [3, \infty) \rightarrow \mathbb{R}$ such that $q(y) = f(x',y)$ for $y \geq 3$. We have $x' \geq x_1 - \frac{1}{2} > 3$ and $q(3) = \left( \frac{5}{6-x'} \right)^{5/2} - x' > \left( \frac{5}{6-3} \right)^{5/2} - x_1 > 0$. For any $y \in [3, \infty)$, let $z_y  = \frac{2y - x'}{x'-1}$. As $y$ increases to infinity, $z_y$ increases to infinity. We have $q(y) = \left( \left(1 + \frac{1}{z_y} \right)^{z_y +1}\right)^{(x'-1)/2} - x'$. Thus, by Lemma~\ref{limitlemma}, $q(y)$ decreases from $q(3) > 0$ to $h(x') < 0$ as $y$ increases from $3$ to infinity. Thus, there exists some $y_{\delta} \in [3, \infty)$ such that $q(y_{\delta}) = 0$. We have $f(x',y_{\delta}) = 0$, so $x' = x_{y_{\delta}}$. By (a), $x_y > x_{y_{\delta}}$ for any $y \in (y_{\delta}, \infty)$. Thus, (b) is proved.~\hfill{$\Box$}
\\
\\
\textbf{Proof of Theorem~\ref{compresulte}.} 
Let $x = kt/m$. Since $m < kt$, $x > 1$. By (\ref{eqnsve5}), $m \geq kt/2$, so $x \leq 2$. Let $\sim$ be any of the relations $<$, $=$, and $>$. We have
\begin{align} &b_1(k,t,m) \sim b_2(k,t,m) \quad \Leftrightarrow \quad \frac{m + (k-1)t}{m} \sim (2k-1) \left ( 1- \frac{k-1}{k} { \Big( \frac{m}{kt} \Big) }^{1/(k-1)} \right ) \nonumber \\
&\Leftrightarrow \quad 1 + (k-1)\frac{x}{k} \sim 2k - 1 - \frac{(2k-1)(k-1)}{k x^{1/(k-1)}} \nonumber \\
&\Leftrightarrow \quad \frac{(2k-1)(k-1)}{x^{1/(k-1)}} \sim 2k(k-1) - (k-1)x \quad \Leftrightarrow \quad \frac{(2k-1)}{x^{1/(k-1)}} \sim 2k - x \nonumber \\
&\Leftrightarrow \quad (2k-1)^{k-1} \sim (2k-x)^{k-1} x > 0 \quad \mbox{(as $1< x \leq 2$)} \nonumber \\
&\Leftrightarrow \quad \left(\frac{2k-1}{2k-x} \right)^{k-1} \sim x. \label{fact2e} 
\end{align} 

Let $f$ be as in Lemma~\ref{limitlemma2e}. Let $y_0 = k$. By Lemma~\ref{limitlemma2e}, $f(x_{y_0}, y_0) = 0$ for some unique $x_{y_0} \in (1, x_1)$, and the larger $y_0$ is, the larger $x_{y_0}$ is. Let $x_0 = x_{y_0}$. It can be checked that $x_0 = 2.088...$ if $k = 3$. 
By Lemma~\ref{limitlemma2e} (b), for any real $\delta > 0$, $x_0 > x_1 - \delta$ if $y_0$ is sufficiently large.  
Suppose $m \geq kt/x_0$. Then, $x \leq x_0$. We have $\left(\frac{2k-1}{2k-x} \right)^{k-1/2} - x = f(x,y_0) \leq 0$ by Lemma~\ref{limitlemma2e}. Since $\left(\frac{2k-1}{2k-x} \right)^{k-1} < \left(\frac{2k-1}{2k-x} \right)^{k-1/2} \leq x$, 
%
%
we have $b_1(k,t,m) < b_2(k,t,m)$ by (\ref{fact2e}). Thus, (a) is proved.

We now prove (b). Let $z = \frac{2k - x}{x-1}$ and $c = \left(\frac{2k-2}{2k-1}\right)^{1/2}$. 
We have 
\begin{align} \left(\frac{2k-1}{2k-x} \right)^{k-1} &= \left(\frac{2k-1}{2k-x} \right)^{k-1/2} \left(\frac{2k-1}{2k-x}\right)^{-1/2} = \left(1 + \frac{1}{z} \right)^{(z+1)(x-1)/2} \left(\frac{2k-x}{2k-1}\right)^{1/2}  \nonumber \\
&\geq c \left(\left(1 + \frac{1}{z} \right)^{(z+1)}\right)^{(x-1)/2}  \quad \mbox{(as $x \leq 2$)}\nonumber \\
&> c e^{(x-1)/2} \quad \mbox{(by Lemma~\ref{limitlemma}).} 
\label{fact3e}
\end{align}
%
Let $f_c$ be as in Lemma~\ref{limitlemma3}. By Lemma~\ref{limitlemma3}, $f_c(x_0') = 1$ for some unique $x_0' \in (x_1,\infty)$, and the larger $c$ is, the smaller $x_0'$ is. Thus, the larger $k$ is, the smaller $x_0'$ is. It can be checked that $x_0' = 3.991...$ if $k = 3$. 
By Lemma~\ref{limitlemma3} (d), for any real $\delta > 0$, $x_0' < x_1 + \delta$ if $k$ is sufficiently large. Suppose $m \leq kt/x_0'$. Then, $x \geq x_0'$. By Lemma~\ref{limitlemma3}, we have $f_c(x) \geq 1$, so $c e^{(x-1)/2} \geq x$. By (\ref{fact2e}) and (\ref{fact3e}), $b_1(k,t,m) > b_2(k,t,m)$.~\hfill{$\Box$}

\section{The bounds in Theorems~\ref{v1} and \ref{v2}} \label{comp1section}

Let $k$, $t$, and $n$ be as in Theorems~\ref{v1} and \ref{v2}. Let $b_1(k,t,n)$ and $b_2(k,t,n)$ be the bound in Theorem~\ref{v1} and the bound in Theorem~\ref{v2}, respectively; that is,
\[b_1(k,t,n) = \frac{n + (k-1) t}{2k} \quad \mbox{and} \quad b_2(k,t,n) = n \left( 1 - \frac{k}{k+1} \left( \frac{n}{(k+1)t} \right)^{1/k} \right).\]
%
If $k = 1$, then clearly $n = t$, $b_1(k,t,n)$ is attained, and $b_1(k,t,n) < b_2(k,t,n)$. We now consider $k \geq 2$. By (\ref{eqnsve3}), $n \leq (k+1) t$, and equality holds if $G$ is the union of $t$ pairwise vertex-disjoint $(k+1)$-vertex stars, in which case the two bounds are equal and attained. We now consider $n < (k+1) t$. 


\begin{theorem} \label{compresult} Suppose $k \geq 2$ and $n < (k+1)t$.\medskip 
\\
(a) Let $x_1 > 1$ be the real number $3.512...$ such that $e^{(x_1 - 1)/2} = x_1$. There exists a unique real number $x_0$ such that $2.438... \leq x_0 < x_1$ and $\left( \frac{2k}{2k + 1 - x_0} \right)^{k} = x_0$. We have 
\[b_1(k,t,n) < b_2(k,t,n) \quad \mbox{if } \; n \geq \frac{1}{x_0} (k+1) t.\]
The larger $k$ is, the larger $x_0$ is. Moreover, for any real $\delta > 0$, $x_0 > x_1 - \delta$ if $k$ is sufficiently large.
\medskip 
\\
(b) Let $x_{1/4} > 1$ be the real number $7.908...$ such that $\frac{1}{4} e^{(x_{1/4} - 1)/2} = x_{1/4}$. There exists a unique real number $x_0'$ such that $5.594... \leq x_0' < x_{1/4}$ and $\left( \frac{k+4}{2k+4} \right)^2 e^{(x_0' - 1)/2} = x_0'$. We have 
\[b_2(k,t,n) < b_1(k,t,n) \quad \mbox{if } \; n \leq \frac{1}{x_0'} (k+1)t.\]
The smaller $k$ is, the smaller $x_0'$ is. 
\end{theorem}
Since $n$ can be at most $(k+1)t$, this result tells us that the range of values of $n$ for which the bound in Theorem~\ref{v1} is better than the bound in Theorem~\ref{v2} is wider than that for which the opposite holds. 

We now prove Proposition~\ref{compresult}. By slightly modifying the function $f$ in Lemma~\ref{limitlemma2e}, we obtain the following lemma by the same argument for Lemma~\ref{limitlemma2e}. 

\begin{lemma} \label{limitlemma2} Let $A = \{(x, y) \in \mathbb{R} \times \mathbb{R} \colon y \geq 2, \, 1 \leq x < 2y+1\}$. Let $f: A \rightarrow \mathbb{R}$ be the function given by
\[f(x,y) = \left( \frac{2y}{2y + 1 - x} \right)^y - x\] 
for $(x,y) \in A$. For any $y_0 \in [2, \infty)$, $f(x_{y_0}, y_0) = 0$ for some unique $x_{y_0} \in (1, 2y_0 + 1)$, and $f(x,y) < 0$ for any $x \in (1,x_{y_0}]$ and $y \in [y_0, \infty)$ such that $x \neq x_{y_0}$ or $y \neq y_0$. 

Moreover, let $x_1 > 1$ be the real number $3.512...$ such that $e^{(x_1 - 1)/2} = x_1$. \\
(a) If $y_0, y_1 \in [2, \infty)$ with $y_0 < y_1$, then $x_{y_0} < x_{y_1} < x_1$. \\
(b) For any real $\delta > 0$, there exists some $y_{\delta} \in [2, \infty)$ such that $x_y > x_1 - \delta$ for any $y \in (y_{\delta}, \infty)$.
\end{lemma}
\textbf{Proof of Theorem~\ref{compresult}.} 
Let $x = (k+1) t/n$. Since $n < (k+1) t$, $x > 1$. Obviously, $n \geq t \geq 1$, so $x \leq k+1$. Let $\sim$ be any of the relations $<$, $=$, and $>$. We have
\begin{align} &b_1(k,t,n) \sim b_2(k,t,n) \quad \Leftrightarrow \quad \frac{n + (k-1) t}{n} \sim 2k \left( 1 - \frac{k}{k+1} \left( \frac{n}{(k+1)t} \right)^{1/k} \right) \nonumber \\
&\Leftrightarrow \quad 1 + \frac{(k-1)x}{k+1} \sim 2k - \frac{2k^2}{(k+1)x^{1/k}} \quad \Leftrightarrow \quad \frac{2k^2}{x^{1/k}} \sim (2k-1)(k+1) - (k-1)x \nonumber \\
&\Leftrightarrow \quad (2k^2)^{k} \sim ((2k-1)(k+1) - (k-1)x)^k x > 0 \quad \mbox{(as $x \leq k+1$)} \nonumber \\
&\Leftrightarrow \quad \left(\frac{2k^2}{2k^2 - (k-1)(x-1)} \right)^{k} \sim x. \label{fact2} 
\end{align} 

Let $f$ be as in Lemma~\ref{limitlemma2}. Let $y_0 = k$. Let $x_1$ be as in (a). By Lemma~\ref{limitlemma2}, $f(x_{y_0}, y_0) = 0$ for some unique $x_{y_0} \in (1, x_1)$, and the larger $y_0$ is, the larger $x_{y_0}$ is. Let $x_0 = x_{y_0}$. It can be checked that $x_0 = 2.438...$ if $k = 2$. 
By Lemma~\ref{limitlemma2} (b), for any real $\delta > 0$, $x_0 > x_1 - \delta$ if $y_0$ is sufficiently large.  
Suppose $n \geq (k+1) t/x_0$. Then, $x \leq x_0$. We have $\left(\frac{2k}{2k + 1 - x} \right)^{k} - x = f(x,y_0) \leq 0$ by Lemma~\ref{limitlemma2}. Since 
\[\left( \frac{2k^2}{2k^2 - (k-1)(x-1)} \right)^k < \left( \frac{2k^2}{2k^2 - k(x-1)} \right)^k = \left(\frac{2k}{2k + 1 - x} \right)^{k} \leq x,\]
we have $b_1(k,t,n) < b_2(k,t,n)$ by (\ref{fact2}). Thus, (a) is proved.

We now prove (b). Let $z = \frac{2k + 5 - x}{x-1}$ and $c = \left( \frac{k + 4}{2k+4} \right)^2$. Since $x \leq k+1$, $c \leq \left(\frac{2k + 5 - x}{2k+4} \right)^2$. We have 
\begin{align} &\left(\frac{2k^2}{2k^2 - (k-1)(x-1)} \right)^{k} \geq \left(\frac{2(k-1)(k+2)}{2(k-1)(k+2) - (k-1)(x-1)} \right)^{k}\nonumber \\
&= \left(\frac{2(k+2)}{2(k+2) - (x-1)} \right)^{k} = \left(\frac{2k+4}{2k + 5 - x} \right)^{k+2} \left(\frac{2k + 5 - x}{2k+4} \right)^{2} \nonumber \\
&\geq c \left(\frac{2k+4}{2k + 5 - x} \right)^{k+2} = c \left(1 + \frac{x - 1}{2k + 5 - x} \right)^{(z+1)(x-1)/2} = c \left(\left(1 + \frac{1}{z} \right)^{z+1} \right)^{(x-1)/2} \nonumber \\
&> c e^{(x-1)/2} \quad \mbox{(by Lemma~\ref{limitlemma}).} 
\label{fact3}
\end{align}
%
Let $f_c$ be as in Lemma~\ref{limitlemma3}. We have $c > \left( \frac{1}{2} \right)^2 = \frac{1}{4}$. By Lemma~\ref{limitlemma3}, $f_c(x_0') = 1$ for some unique $x_0' \in (1, x_{1/4})$, and the larger $c$ is, the smaller $x_0'$ is. Thus, the smaller $k$ is, the smaller $x_0'$ is. It can be checked that $x_0' = 5.594...$ if $k = 2$. 
Suppose $n \leq (k+1)t/x_0'$. Then, $x \geq x_0'$. By Lemma~\ref{limitlemma3}, we have $f_c(x) \geq 1$, so $c e^{(x-1)/2} \geq x$. By (\ref{fact2}) and (\ref{fact3}), $b_1(k,t,n) > b_2(k,t,n)$.~\hfill{$\Box$}

\end{document}